\documentclass[a4paper,11pt]{article}

\usepackage{amssymb}

\usepackage[english,francais]{babel}

\usepackage{amsmath,amsfonts,amssymb,mathtools,stmaryrd,mathrsfs}
\usepackage{dsfont}

\newtheorem{theorem}{Theorem}[section]

\newtheorem{e-proposition}[theorem]{Proposition}
\newtheorem{corollary}[theorem]{Corollary}
\newtheorem{e-definition}[theorem]{Definition\rm}

\newtheorem{theoreme}{Th\'eor\`eme}[section]

\newtheorem{proposition}[theoreme]{Proposition}

\newcommand{\tore}{\mathbb{T}}
\newcommand{\R}{\mathbb{R}}
\newcommand{\N}{\mathbb{N}}
\newcommand{\Z}{\mathbb{Z}}

\newcommand{\ph}{\varphi}
\newcommand{\md}{m}

\usepackage{authblk}
\usepackage{hyperref}

\author[1]{Michel Bena\"im \thanks{michel.benaim@unine.ch,brehier@math.univ-lyon1.fr}}
\author[2]{Charles-Edouard
  Br\'ehier}
\affil[1]{Universit\'e de Neuch\^atel,  Institut de Math\'ematiques, Rue Emile Argand 11, CH-2000 Neuch\^atel, Switzerland}
\affil[2]{Univ Lyon, Universit\'e Claude Bernard Lyon 1, CNRS UMR 5208, Institut Camille Jordan, 43 blvd. du 11 novembre 1918, F-69622 Villeurbanne cedex, France}

\title{Convergence of Adaptive Biasing Potential methods for diffusions}

\begin{document}

\maketitle

\begin{abstract}
We prove the consistency of an adaptive importance sampling strategy based on biasing the potential energy function $V$ of a diffusion process $dX_t^0=-\nabla V(X_t^0)dt+dW_t$; for the sake of simplicity, periodic boundary conditions are assumed, so that $X_t^0$ lives on the flat $d$-dimensional torus. The goal is to sample its invariant distribution $\mu=Z^{-1}\exp\bigl(-V(x)\bigr)\,dx$. The bias $V_t-V$, where $V_t$ is the new (random and time-dependent) potential function, acts only on some coordinates of the system, and is designed to flatten the corresponding empirical occupation measure of the diffusion $X$ in the large time regime.

The diffusion process writes $dX_t=-\nabla V_t(X_t)dt+dW_t$, where the bias $V_t-V$ is function of the key quantity $\overline{\mu}_t$: a probability occupation measure which depends on the past of the process, {\it i.e.} on $(X_s)_{s\in [0,t]}$. We are thus dealing with a self-interacting diffusion.

In this note, we prove that when $t$ goes to infinity, $\overline{\mu}_t$ almost surely converges to $\mu$. Moreover, the approach is justified by the convergence of the bias to a limit which has an intepretation in terms of a free energy.

The main argument is a change of variables, which formally validates the consistency of the approach. The convergence is then rigorously proven adapting the ODE method from stochastic approximation.
\end{abstract}

\section{Introduction}
\label{}

Computing the average $\mu(\ph)=\int_{\mathcal{D}}\ph(x)\mu(dx)$ of a function $\ph:\mathcal{D}\to \R$, with respect to a probability distribution $\mu$ defined on $\mathcal{D}\subset \R^d$, is typically a challenging task in many applications (e.g. chemistry, statistical physics, see e.g.~\cite{LeimkuhlerMatthews}), since usually $d$ is large and $\mu$ is multimodal.

In the sequel, we assume that $\mathcal{D}=\tore^d=(\R/\Z)^d$ is the flat $d$-dimensional torus, and that $\mu$ writes
\begin{equation}\label{eq:mu_simple}
\mu(dx)=\mu_\beta(dx)=\frac{\exp\bigl(-\beta V(x)\bigr)}{Z(\beta)}dx,
\end{equation}
where $V:\tore^d\to \R$ is a smooth potential function, $\beta\in(0,+\infty)$ is the inverse temperature, $dx$ denotes the Lebesgue measure on $\tore^d$ and $Z(\beta)$ is a normalizing constant. 
In this context, the multimodality of $\mu_{\beta}$ follows, in the case of so-called energetic barriers, from the existence of several local minima of $V$.


A standard approach to computing $\mu_\beta(\ph)$ is to consider the following SDE on $\tore^d$ (overdamped Langevin dynamics):
\begin{equation}\label{eq:OL_simple}
dX_t^0=-\nabla V(X_t^0)dt+\sqrt{2\beta^{-1}}dW_t, \quad X_0^0=x.
\end{equation}
where $\bigl(W(t)\bigr)_{t\ge 0}$ is standard Brownian Motion on $\tore^d$. Indeed, it is well-known that, for any continuous function $\ph:\tore^d\to \R$ almost surely
\begin{equation}\label{eq:ergo_simple}
\frac{1}{t}\int_{0}^{t}\ph(X_r^0)dr\underset{t\to +\infty}\to \int_{\tore^d}\ph(x)\mu_\beta(dx),
\end{equation}
However, this convergence may be very slow, when $\beta$ is large and $V$ has several minima: the stochastic process $X^0$ is then metastable, and hopping from the neighborhood of one local minimum of $V$ to another is a rare event which may have a strong influence on the estimation of averages $\mu_\beta(\ph)$.

Many strategies based on importance sampling techniques --  self-healing umbrella-sampling \cite{marsili2006self}, well-tempered metadynamics \cite{barducci2008well}, Wang-Landau algorithms, adaptive biasing force, etc... -- have been proposed and applied to improve the convergence to equilibrium of stochastic processes in order to compute approximations of $\mu_{\beta}$. We refer for instance to~\cite{LelievreRoussetStoltz:10} and references therein for a mathematical review.

In this work, we focus on an Adaptive Biasing Potential (ABP) method, given by the system~\eqref{eq:ABP_simple}. The method was designed in~\cite{DicksonLegollLelievreStoltzFleurat:10,marsili2006self} for problems in chemistry, and up to our knowledge no rigorous general mathematical analysis has been performed so far. Precisely, in~\eqref{eq:OL_simple}, $V$ is replaced with a time-dependent and random potential function $V_t$ which is modified adaptively, using the history of the process up to time~$t$: $A_t$ depends on the values of the associated stochastic process $X_r$ for all $0\le r\le t$. Here, $V_t=V-A_t\circ\xi$, where, for some $\md\in\left\{1,\ldots,d-1\right\}$, $A_t:\tore^m\to \R$ and $\xi:\tore^d\to \tore^{\md}$ is a smooth function,  referred to as the {\it reaction coordinate} mapping. In applications, usually $\md\in\left\{1,2,3\right\}$. To simplify further the presentation, we assume that $\xi(x_1,\ldots,x_d)=(x_1,\ldots,x_{\md})$; in this case, $z=(x_1,\ldots,x_{\md})=\xi(x_1,\ldots,x_d)$ (resp. $z^\perp=(x_{\md+1},\ldots,x_d)$) is interpreted as the slow (resp. fast) variable.

The dynamics of the ABP method is given by the following system
\begin{equation}\label{eq:ABP_simple}
\begin{cases}
dX_t=-\nabla\bigl(V-A_t\circ \xi\bigr)(X_t)dt+\sqrt{2\beta^{-1}}dW(t)\\
\overline{\mu}_t=\frac{\overline{\mu}_0+\int_{0}^{t}\exp\bigl(-\beta A_r\circ\xi(X_r)\bigr)\delta_{X_r}dr}{1+\int_{0}^{t}\exp\bigl(-\beta A_r\circ\xi(X_r)\bigr)dr}\\
\exp\bigl(-\beta A_t(z)\bigr)=\int_{\tore^d}K\bigl(z,\xi(x)\bigr)\overline{\mu}_t(dx),~\forall z\in\tore^m,
\end{cases}
\end{equation}
where a smooth {\it kernel} function $K:\tore^\md\times \tore^\md \to (0,+\infty)$, which is such that $\int_{\tore^\md}K(z,\zeta)dz=1, \forall \zeta\in \tore^\md$, is introduced to define a smooth function $A_t$ from the distribution $\overline{\mu}_t$. The unknows in~\eqref{eq:ABP_simple} are the stochastic processes $t\mapsto X_t \in\tore^d$, $t\mapsto \overline{\mu}_t\in \mathcal{P}(\tore^d)$ (the set of Borel probability distributions on $\tore^d$, endowed with the usual topology of weak convergence of probability distributions), and $t\mapsto A_t\in \mathcal{C}^{\infty}(\tore^\md)$ (the set of infinitely differentiable functions on $\tore^\md$). In addition to~\eqref{eq:ABP_simple}, arbitrary (and deterministic, for simplicity) initial conditions $X_{t=0}=x$, $\overline{\mu}_{t=0}=\overline{\mu}_0$ and $A_{t=0}=A_0$ are prescribed.

The third equation in~\eqref{eq:ABP_simple} introduces a coupling between the evolutions of the diffusion $X_t$ and of the {\it weighted empirical distribution} $\overline{\mu}_t$: then $X$ can be seen as a self-interacting diffusion process, like in~\cite{BenaimLedouxRaimond1}.

%


Our main result is the consistency of the ABP approach.
\begin{theorem}\label{th:cv_mu}
Almost surely, $\overline{\mu}_t\underset{t\to +\infty}\to\mu_{\beta}$, in $\mathcal{P}(\tore^d)$.
\end{theorem}

With standard arguments, Theorem~\ref{th:cv_mu} yields almost sure convergence of $A_t$ in $\mathcal{C}^{k}(\tore^\md)$, for all $k\in\N$.
\begin{corollary}\label{cor:cv_A}
Set $\exp\bigl(-\beta A_\infty\bigr)=\int K(\cdot,\xi(x))\mu_{\beta}(dx)$. Then almost surely, $A_t\underset{t\to +\infty}\to A_\infty$, in $\mathcal{C}^{k}(\tore^\md)$, $\forall~k\in\N$.
\end{corollary}

The limit $A_\infty$ is an approximation of the function known as the {\it free energy} $A_\star$ (see~\eqref{eq:A_simple}), which depends on $V$, $\beta$ and $\xi$. As explained in Section~\ref{sec:FE}, the construction of the adaptive dynamics~\eqref{eq:ABP_simple} is motivated by an efficient non-adaptive biasing method, \eqref{eq:bias_simple}, which depends on $A_\star$. Computing $A_\star$ is the aim of many algorithms in molecular dynamics (see~\cite{LelievreRoussetStoltz:10}), and adaptive methods are among the most used in practice. Our result, Theorem~\ref{th:cv_mu}, answers positively the important question of the consistency of ABP method.



The remaining part of the article is organized as follows. In Section~\ref{sec:FE}, we define the free energy function $A_\star$, and explain why non-adaptive and adaptive biaising methods which are related to this function are interesting in the context of metastable dynamics~\eqref{eq:OL_simple}. In Section~\ref{sec:ABP}, we detail the strategy for the proof of Theorem~\ref{th:cv_mu}: we prove a stability estimate for $A_t$, and then introduce a random change of variables, based on a change of time. We are then in position to adapt the strategy of proof from~\cite{BenaimLedouxRaimond1} in our setting, which is based on the ODE method from stochastic approximation. The main essential role of the change of variables is the identification of the limit flow.

The main result Theorem~\ref{th:cv_mu} holds in a more general setting, with appropriate modifications, than that of the present paper. For instance, the overdamped Langevin dynamics may be defined on the non-compact space $\R^d$ instead of $\tore^d$; one can also consider the (hypoelliptic) Langevin dynamics, or infinite-dimensional dynamics (parabolic SPDEs). It is also possible to study the efficiency of the method in terms of a Central Limit Theorem. These generalizations will be studied in~\cite{BenaimBrehier}.


\section{Free energy and construction of the ABP dynamics~\eqref{eq:ABP_simple}}\label{sec:FE}

The aims of this section are to explain first how the ABP method~\eqref{eq:ABP_simple}, is constructed in a consistent way (the limit in Theorem~\ref{th:cv_mu} is $\mu_\beta$); and second why it is expected to be efficient (a rigorous analysis of the efficiency is out of the scope of this work).

Observe that $\exp\bigl(-\beta A_\infty(z)\bigr)=\int_{\tore^m}K(z,\zeta)\exp\bigl(-\beta A_{\star}(\zeta,\beta)\bigr)d\zeta$, where $A_\star(\cdot,\beta)$ is the free energy (at temperature $\beta^{-1}$), defined by: for all $z\in \tore^\md$
\begin{equation}\label{eq:A_simple}
\exp\bigl(-\beta A_{\star}(z,\beta)\bigr)=\int_{\tore^{d-\md}}\frac{\exp\bigl(-\beta V(z,z^\perp)\bigr)}{Z(\beta)}dz^\perp.
\end{equation}
Usually, $K(z,\zeta)=K^\epsilon(z,\zeta)=\frac{1}{\epsilon}\tilde{K}\bigl((\zeta-z)/\epsilon\bigr)$, where $\epsilon\in(0,1)$ and $\tilde{K}:\R^\md\to (0,+\infty)$ is symmetric, smooth, with compact support in $[-1/2,1/2]$; then $A_\infty^\epsilon$ converges to $A_\star(\cdot,\beta)$, in $\mathcal{C}^\infty$. Choosing $\epsilon$ sufficiently small, $A_t$ almost surely approximates the free energy $A_{\star}(\cdot,\beta)$ when  $t\to +\infty$, thanks to Corollary~\ref{cor:cv_A}.

Equation~\eqref{eq:A_simple} means that $\exp\bigl(-\beta A_{\star}(z,\beta)\bigr)dz\in \mathcal{P}(\tore^\md)$ is the image $\mu_\beta\bigl(\xi^{-1}(\cdot)\bigr)$ of $\mu_\beta$ by $\xi$. The free energy gives an effective potential along $\xi$, which is chosen in practice such that $\bigl(\xi(X_t^0)\bigr)_{t\ge 0}$ is metastable; this is related to $\mu_\beta\bigl(\xi^{-1}(\cdot)\bigr)$ being metastable, for instance when $A_{\star}(\cdot,\beta)$ has several local minima.


This is why in many applications, computing free energy differences $A_{\star}(z_1,\beta)-A_\star(z_2,\beta)$ is essential, see~\cite{LelievreRoussetStoltz:10}. The free energy function also theoretically provides efficient importance sampling algorithms; however these algorithms can only be implemented if $A_\star$ is explicitly known, and adaptive strategies allow to circumvent this practical difficulty.  Define biased probability distribution and dynamics
\begin{equation}\label{eq:bias_simple}
\begin{gathered}
\mu_{\beta}^{{\star}}=\frac{\exp\bigl(-\beta\bigl[V(x)-A_{\star}(\xi(x),\beta)\bigr]\bigr)}{Z(\beta)}dx\\
dX_{t}^{\star}=-\nabla\bigl[V-A_{\star}(\xi(\cdot),\beta)\bigr](X_{t}^{\star})dt+\sqrt{2\beta^{-1}}dW(t),
\end{gathered}
\end{equation}
by replacing the original potential function $V$ with the biased potential function $V-A_{\star}\bigl(\xi(\cdot),\beta\bigr)$ in~\eqref{eq:mu_simple} and~\eqref{eq:OL_simple}. Note that $\mu_{\beta}^{\star}$ is the unique invariant distribution of $X^\star$. By construction, it is easy to check that the image by $\xi$ of $\mu_{\beta}^{\star}$ is the uniform distribution $dz$ on $\tore^\md$, {\it i.e.} the associated free energy is equal to $0$.

Now define (unweighted) empirical distributions associated with~\eqref{eq:OL_simple} and~\eqref{eq:bias_simple}  respectively:
$$\rho_t^0=\frac{1}{t}\int_{0}^{t}\delta_{X_r^0}dr \quad,\quad \rho_t^\star=\frac{1}{t}\int_{0}^{t}\delta_{X_r^\star}dr.$$
Then, by~\eqref{eq:ergo_simple}, the image by $\xi$ of $\rho_t^0$, resp. $\rho_t^{\star}$, converges almost surely in $\mathcal{P}(\tore^m)$, to $\exp\bigl(-\beta A_\star(z,\beta)\bigr)dz$, resp. $dz$. Thus the dynamics in~\eqref{eq:bias_simple} reaches asymptotically a {\it flat histogram} property in the $z=\xi(x)$ direction; the exploration of $\tore^\md$ is thus faster for $\xi(X^\star)$ than for $\xi(X^0)$, and in turn the convergence of $X^\star$ to $\mu_{\beta}^{\star}$ is expected to be faster than the convergence of $X^0$ to $\mu_\beta$.

Finally, the construction of the ABP method~\eqref{eq:ABP_simple}, in particular the use of
{\it weighted} empirical distributions $\overline{\mu}_t$, is motivated by the following almost sure convergence: for any continuous $\ph:\tore^d\to \R$,
\begin{equation}\label{eq:identity_bias-unbiased_simple}
\frac{\frac{1}{t}\int_{0}^{t}\exp\bigl(-\beta A_\star(\xi(X_r^\star),\beta)\bigr)\ph(X_r^\star)dr}{\frac{1}{t}\int_{0}^{t}\exp\bigl(-\beta A_\star(\xi(X_r^\star),\beta)\bigr)dr}\underset{t\to +\infty}\to\mu_{\beta}^{\star}\bigl(\ph \exp\bigl(-\beta A_{\star}(\xi(\cdot),\beta)\bigr)\bigr)=\mu_{\beta}(\ph).
\end{equation}
Theorem~\ref{th:cv_mu} thus extends this consistency property from a non-adaptive \eqref{eq:bias_simple} to an adaptive dynamics~\eqref{eq:ABP_simple}.



\section{Proof of Theorem~\ref{th:cv_mu}}\label{sec:ABP}

In this section, we provide the main ideas of the proof of Theorem~\ref{th:cv_mu}.Some technical arguments are skipped, and will be fully detailed in~\cite{BenaimBrehier}, in a more general framework. We first state an important property of $A_t$, and then introduce a change of variables, which helps us identifying a more standard form for self-interacting diffusion processes. We then adapt in our context the arguments from~\cite{BenaimLedouxRaimond1}, to establish the consistency of the ABP approach thanks to the ODE method from stochastic approximation theory.

\subsection{Properties of the ABP dynamics~\eqref{eq:ABP_simple}}

Our first task in the study of the ABP dynamics is to study the well-posedness of the equation, {\it i.e.} the existence of a unique global solution $t\mapsto (X_t,\overline{\mu}_t,A_t)\in \tore^d\times \mathcal{P}(\tore^d)\times \mathcal{C}^\infty(\tore^\md)$. In order to apply a standard fixed point/Picard iteration strategy, it is essential to control the Lipschitz constant of $\nabla\bigl(A_t\circ\xi\bigr)$ (first equation in~\eqref{eq:ABP_simple}). This key stability property is ensured as follows. Let $m=\min_{z,\zeta\in\tore^\md}K(z,\zeta)$, and $M^{(n)}=\max_{z,\zeta}|\partial_z^{n}K(z,\zeta)|$ for $n\in\left\{0,1\right\}$, where $\partial_z^n$ denotes the differential of order $n$, and introduce
$$\mathcal{A}=\left\{A\in \mathcal{C}^{\infty}(\tore^\md)~|~ \min_{z\in\tore^\md}e^{-\beta A(z)}\ge m, \max_{z\in\tore^\md}|\partial_z^n e^{-\beta A(z)}|\le M^{(n)}, n=0,1\right\}.$$ 
Then $\mathcal{A}$ is  left invariant by the evolution $t\mapsto A_t$, {\it i.e.} $A_0\in\mathcal{A}$ implies $A_t\in\mathcal{A}$ for all $t\ge 0$, almost surely. 

%

\subsection{Change of variables}

The stochastic process $t\mapsto \overline{\mu}_t$, with values in $\mathcal{P}(\tore^d)$, is the unique solution of the random Ordinary Differential Equation (ODE), interpreted in a weak sense (considering continuous bounded test functions):
\begin{equation}\label{eq:ODE_mu}
\frac{d\overline{\mu}_t}{dt}=\frac{\theta'(t)}{1+\theta(t)}\bigl(\delta_{X_t}-\overline{\mu}_t\bigr),
\end{equation}
where $\theta(t)=\int_{0}^{t}\exp\bigl(-\beta A_r(\xi(X_r))\bigr)dr$. 
The random function $\theta:[0,+\infty)\to [0,+\infty)$ is a $\mathcal{C}^1$-diffeomorphism: indeed for all $t\ge 0$, almost surely $\theta'(t)=\exp\bigl(-\beta A_t(\xi(X_t))\bigr)\in[m,M]$. This fundamental property allows us to apply the following change of variables:
\begin{equation}\label{eq:change}
s=\theta(t) \quad,\quad t=\theta^{-1}(s) \quad;\quad Y_s=X_t \quad,\quad \overline{\nu}_s=\overline{\mu}_t \quad,\quad B_s=A_t.
\end{equation}
Observe that $s=\theta(t)\underset{t\to +\infty}\to +\infty$ and that $t=\theta^{-1}(s)\underset{s\to +\infty}\to +\infty$, almost surely. Instead of studying the asymptotic behavior of $\overline{\mu}_t$ when $t\to +\infty$, it thus equivalent to study the asymptotic behavior of $\overline{\nu}_s$ when $s\to +\infty$. In the new variables~\eqref{eq:change}, the ABP dynamics~\eqref{eq:ABP_simple} writes
\begin{equation}\label{eq:ABP_simple_s}
\begin{cases}
dY_s=-\nabla\bigl(V-B_s\circ \xi\bigr)(Y_s)e^{\beta B_s(\xi(Y_s))}ds+\sqrt{2\beta^{-1}e^{\beta B_s(\xi(Y_s))}}d\tilde{W}(s)\\
\overline{\nu}_s=\frac{\overline{\nu}_0+\int_{0}^{s}\delta_{Y_r}dr}{1+s}\\
\exp\bigl(-\beta B_s(z)\bigr)=\int_{\tore^d}K(z,\xi(x))\overline{\nu}_s(dx),
\end{cases}
\end{equation}
where $\tilde{W}$ is a new standard Brownian motion on $\tore^d$, defined from $W$ and $\theta$. Notice that $\overline{\nu}_s$ is a {\it nonweighted empirical distribution} and that $s\mapsto \overline{\nu}_s$ satisfies the simpler random ODE
\begin{equation}\label{eq:ODE_nu}
\frac{d\overline{\nu}_s}{ds}=\frac{1}{1+s}\bigl(\delta_{Y_s}-\overline{\nu}_s\bigr).
\end{equation}
The change of variable~\eqref{eq:change} both removes $\theta(t)$ from~\eqref{eq:ODE_mu} as well as the weigths $\exp\bigl(-\beta A_t(\xi(X_t))\bigr)=\theta'(t)$ from~\eqref{eq:ABP_simple}.

Thanks to Equation~\eqref{eq:ABP_simple_s}, an analogy with the framework of~\cite{BenaimLedouxRaimond1} can now be made. Even though we cannot directly apply the results therein, due to the specific form of the dynamics on $Y$, we follow the same strategy for the analysis of $\overline{\nu}_s$ when $s\to +\infty$: we use the {\it ODE method}.

\subsection{Application of the ODE method and sketch of proof of Theorem~\ref{th:cv_mu}}

The guideline of the so-called ODE approach we wish to apply is as follows: there is an asymptotic time-scale separation between the (fast) evolution of $Y_s$ and the (slow) evolution of $\overline{\nu}_s$ (and of $B_s$). The asymptotic behavior of $\overline{\nu}_s$ is then determined by a so-called limit ODE, where $\delta_{Y_s}$ is replaced in~\eqref{eq:ODE_nu} with the unique invariant probability distribution of the following SDE on $\tore^d$,
\begin{equation}\label{eq:Y_frozen_simple}
dY_s^B=-\nabla\bigl(V-B\circ \xi\bigr)(Y_s)e^{\beta B(\xi(Y_s))}ds+\sqrt{2\beta^{-1}e^{\beta B(\xi(Y_s))}}d\tilde{W}(s),
\end{equation}
{\it i.e.} the first (fast) equation of~\eqref{eq:ABP_simple_s} where the slowly varying variable $B_s$ is frozen at an arbitrary $B\in\mathcal{A}$. In fact, we have the following fundamental result: the invariant distribution of~\eqref{eq:Y_frozen_simple} does not depend on $B$.
\begin{proposition}\label{propo:consis}
For any smooth $B:\tore^\md\to \R$, the unique invariant distribution of~\eqref{eq:Y_frozen_simple} is $\mu_\beta$.
\end{proposition}
Proposition~\ref{propo:consis} is essential and its proof is very simple. Indeed, introduce the generator $\mathcal{L}_X^B$ of $X^B$, resp. the unique invariant distribution of $X^B$, denoted by $\mu_\beta^B(dx)=Z^B(\beta)^{-1}\exp\bigl(-\beta (V-B\circ \xi)(x)\bigr)dx$, with the diffusion 
$dX_t^B=-\nabla\bigl(V-B\circ \xi\bigr)(Y_s)dt+\sqrt{2\beta^{-1}}dW(t)$. 
Then the generator $\mathcal{L}_Y^B$ of $Y^B$ defined by~\eqref{eq:Y_frozen_simple} is equal to $\exp\bigl(\beta B\circ\xi\bigr)\mathcal{L}_X^B$. Proposition~\ref{propo:consis} is a consequence of the following identity: for any smooth $\phi,\psi:\tore^d\to \R$,
$$\int_{\tore^d}\phi(y)\mathcal{L}_Y^B\psi(y)\mu_\beta(dy)
=\int_{\tore^d}\phi(x)\mathcal{L}_X^B\psi(x)\mu_\beta^B(dx)=0.$$


We now outline the end of the proof of Theorem~\ref{th:cv_mu}, adapting the arguments from~\cite{BenaimLedouxRaimond1} in our original case; details in a more general setting are given in~\cite{BenaimBrehier}. The ODE method suggests us to define $\Gamma(\sigma,s,\nu)=\Gamma_{\sigma-s}(\nu)$, for any $\sigma\ge s$ and $\nu\in\mathcal{P}(\tore^d)$, where $\Gamma_s(\nu)=e^{-s}\nu+(1-e^{-s})\mu_\beta\underset{s\to +\infty}\to\mu_\beta$ is the solution of $\frac{d\Gamma_s}{ds}=\mu_\beta-\Gamma_s$ with $\Gamma_0(\nu)=\nu$. 
To state (without proof) our last techincal result, we recall that weak convergence in $\mathcal{P}(\tore^d)$ is associated with the following metric
$$d\bigl(\mu^1,\mu^2\bigr)=\sum_{n=1}^{+\infty}\frac{1}{2^n}\min\bigl(1,|\int_{\tore^d}f_nd\mu^1-f_nd\mu^2|\bigr),$$
for a given family $\bigl(f_n\bigr)_{n\ge 1}$ of $\mathcal{C}^\infty$ functions, which is dense in $\mathcal{C}^0(\tore^d)$.

\begin{proposition}\label{propo:APT}
For any $S\ge 0$, almost surely
$$\Delta(s,S)=\sup_{\sigma\in[0,S]} d\bigl(\overline{\nu}_{\exp(s+\sigma)},\Gamma(\sigma,s,\overline{\nu}_s)\bigr)\underset{s\to +\infty}\to 0,$$
{\it i.e.} almost surely $s\mapsto \overline{\nu}_s$ is an {\it asymptotic pseudo-trajectory} of the semi-flow~$\Gamma$.
\end{proposition}
We refer to~\cite{BenaimLedouxRaimond1} for a proof of a similar result in a different context, and to~\cite{BenaimBrehier} for a detailed proof in a more general context; the main difference between the two situations is the use of a specific Poisson equation related to the generator of~\eqref{eq:Y_frozen_simple}.

To conclude, observe that $d\bigl(\overline{\nu}_{\exp(s)},\mu_{\beta}\bigr)\leq \Delta(s-S,S)+d\bigl(\Gamma_{S}(\overline{\nu}_{\exp(s)}),\mu_\beta\bigr)$. Letting first $s$, then $S$, go to $+\infty$, Proposition~\ref{propo:APT} implies the main result of this paper, Theorem~\ref{th:cv_mu}.




\section*{Acknowledgements}

The authors would like to thank Tony Leli\`evre and Gabriel Stoltz for helpful comments. The work is partially supported by the Swiss National Foundation, Grants: $200020 \_ 149871$ and $200021 \_ 163072$.

%

\end{document}